\documentclass[12pt]{elsart}
\usepackage{indentfirst,latexsym}
\usepackage{flafter}
\usepackage{graphicx}
\usepackage{amssymb}
\begin{document}
\runauthor{Shu-Zhen Lai,Hou-Biao Li etc.}
\begin{frontmatter}
\title{A Symmetric Rank-one Quasi Newton Method for Non-negative Matrix Factorization\thanksref{b}}
\thanks[b]{\small Supported by National Natural Science Foundation of China (11026085, 11101071, 1117105, 11271001, 51175443) and the Fundamental Research Funds for China Scholarship Council.}

\author[huang]{Shu-Zhen Lai\corauthref{a}},
\corauth[a]{Corresponding author.}\ead{lihoubiao0189@163.com or
laishuzhen08@126.com}
\author[huang]{Hou-Biao Li},
\author[zhang]{Zu-Tao Zhang},
\address[huang]{School of  Mathematics Science, University of Electronic Science and Technology of
China, Chengdu,  611731, P. R. China}
\address[zhang]{School of Mechanical Engineering, Southwest Jiaotong University,
Chengdu 610031, P. R. China }

\begin{abstract}
As we all known, the nonnegative matrix factorization (NMF) is a
dimension reduction method that has been widely used in image
processing, text compressing and signal processing etc. In this
paper, an algorithm for nonnegative matrix approximation is
proposed. This method mainly bases on the active set  and the
quasi-Newton type algorithm, by using the symmetric rank-one and
negative curvature direction technologies to approximate the Hessian
matrix. Our method improves the recent results of those methods in
[Pattern Recognition, 45(2012) 3557-3565; SIAM J. Sci. Comput.,
33(6)(2011)3261-3281; Neural Computation, 19(10)(2007)2756-2779,
etc]. Moreover, the object function decreases faster than many other
NMF methods. In addition, some numerical experiments are presented
in the synthetic data, imaging processing and text clustering. By
comparing with the other six nonnegative matrix approximation
methods, our experiments confirm to our analysis.
\end{abstract}


\end{frontmatter}

\setcounter{equation}{0}
\renewcommand{\theequation}{1.\arabic{equation}}

\section{Introduction}

An NMF problem is to decompose a nonnegative matrix $V\in\mathbb{R}^{n\times m}$ into two nonnegative matrix $W\in\mathbb{R}^{n\times k}$ and $H\in\mathbb{R}^{k\times m}$, such that the $WH$ approximate to $V$ as well as possible. To measure the distance between $WH$ and $V$, there are many methods  such as the Kullback Leibler divergence, Bregman divergence, Frobenius divergence, etc. Due to the
favorable property of the Frobinius divergence, many methods are presented based on it.  The Frobenius divergence is as follows:
\begin{equation}\label{Eq:1.1}
\min_{(W\in\mathbb{R}^{(n,k)},H\in\mathbb{R}^{(k,m)})}f(W,H)=\frac{1}{2}||V-WH||^2_F
\end{equation}
subject to $W_{ij}\geq 0$ ,$H_{ij}\geq0$ for all $i$ and $j$.

In the last decade, numerous methods have been proposed to deal with the NMF problem in Eq.(\ref{Eq:1.1}). Most of these methods can be classified into two classes i.e., alternating one-step gradient descent and alternating least squares.

\begin{itemize}
  \item The alternating one-step gradient is to alternating $W$ and $H$ with one step. The most well known one is Lee's multiplicative update algorithm \cite{Lee2}, which alternates $W$ and $H$ by the following rules: Suppose we have obtained the $l$th matrix $W^l$ and $H^l$, then
\begin{equation}
W_{ia}^{l+1}=W_{ia}^l\frac{(V(H^l)^T)_{ia}}{(W^lH^l(H^l)^T)_{ia}},\forall i,a,
\end{equation}
\begin{equation}
H_{bj}^{l+1}=H_{bj}^l\frac{((W^{l+1})^TV)_bj}{((W^{l+1})^TW^{l+1}H^k)_{bj}},\forall b,j.
\end{equation}
  \item The frame work of the alternating least squares \cite{Berry} is as follows:

  (1)  Initialize $H\in\mathbb{R}^{k\times m}$ with a nonnegative matrix;\\
\indent (2)  Solving the following problem repeatedly until a stopping criterion is satisfied:
\begin{equation}\label{Eq:1.4}
\min_{W\geq 0}g(W)=\frac{1}{2}||H^TW^T-V^T||_F^2,
\end{equation}
  where $H$ is fixed,and
\begin{equation}\label{Eq:1.5}
\min_{H\geq 0}g(H)=\frac{1}{2}||V-WH||_F^2,
\end{equation}
  where $W$ is fixed.
\end{itemize}

Obviously, the above two methods both satisfied
\begin{equation}
f(W^{l+1},H^l)\leq f(W^l,H^l),
\end{equation}
\begin{equation}
f(W^{l+1},H^{l+1})\leq f(W^{l+1},H^l).
\end{equation}

Since Eq.(\ref{Eq:1.4}) and Eq.(\ref{Eq:1.5}) are regarded as the
subproblems of the NMF. Next we focus on Eq.(\ref{Eq:1.5}),while the
method for Eq.(\ref{Eq:1.4}) is the same as the method of
Eq.(\ref{Eq:1.5}), due to the symmetric property of Euclid distance.
We have learnt some method to solve Eq.(\ref{Eq:1.5}), such as the
projected Newton method \cite{Lin}, quasi-Newton method
\cite{Zdunk,Kim1,Kim2}, active set method \cite{HKim,Jkim1,Jkim2}.
Most of the methods are based on the solution of nonnegative least
squares.

Note that the Frobenius norm of a matrix is just the sum of the Euclidiean distance over columns (or rows). Then solving Eq(\ref{Eq:1.5}) can be boiled down to solve a series of Nonnegative least squares (NNLS) problems of the following forms:
\begin{equation}\label{nnlses}
{\begin{array}{ll}\min_{h_1\geq 0}g(h_1)=\frac{1}{2}||Wh_1-V_1||_2^2,\min_{h_2\geq0}g(h_2)=\frac{1}{2}||Wh_2-V_2||_2^2,\cdots,\\
\min_{h_k\geq 0}g(h_k)=\frac{1}{2}||Wh_k-V_k||_2^2,\end{array}}
\end{equation}
where $ h_i,V_i,i=1,2,\cdots,m$ are the columns of $H$ and $V$, respectively.

Next, for convenience, we write the least squares by omitting their subscripts. For example,
 \begin{equation}\label{Eq:1.8}
 \min g(h)=\frac{1}{2}||Wh-v||_2^2,
 \end{equation}
 subject to $h\geq 0$ for all $i$.

 If $x$ is the optimal solution of Eq.(\ref{Eq:1.8}), then $x$ satisfies the following equation:
 \begin{equation}\label{KKT}
 \left\{{\begin{array}{ll}
 r=W^TWh-W^Tv,\\
 r\geq 0,\\
 x\geq 0,\\
 x^Tv=0.
 \end{array}}\right.
 \end{equation}
 which is called the KKT optimal condition \cite{Bertsekas2}.

 This tells us that we can try to find a decent method which after finite iterations the numerical solution will satisfy
 the condition (\ref{KKT}).

 The remained part of this paper is organized as follows. In the second part we will propose our algorithm for
 a single right hand vector non-negative least squares. In addition, the symmetric rank-one quasi-Newton method for
 NMF will be discussed in the third part. In the fourth part we present the numerical experiments on both synthetic
 data and real world data. And the last part is to be the conclusion part, in this part we will forecast the
 future work for the NMF.
 \setcounter{equation}{0}
\renewcommand{\theequation}{2.\arabic{equation}}

\section{A New Method for Nonnegative Least Squares Problems}

The non-negative least squares can be regarded as a quadratic
programming with box constrains, whose supper constrains can be
considered infinite. Recently, C. Lin proposed the projected
algorithm for NNLS in \cite{Lin}. Later in \cite{Kim1}, D. Kim, S.
Sra and I. S. Ddillon have proposed the fnmae method for
non-negative least square programming based on the quasi-Newton  and
active set method. Few years later, PingHua Gong \cite{gong}
improved the fnmae algorithm by using the symmetric property of the
Hessian matrix of the object function $f(W,H)$. In this paper we
will use the symmetric rank-one quasi-Newton line search technology
approximate the Hessian matrix. In addition, we will combine the
block active method with the quasi Newton method, then use the
symmetric rank-one technology to modify the BFGS in Lin \cite{Lin}.

\subsection{A Symmetric Rank-one Quasi Newton Method}

For the NNLS problem Eq.(\ref{Eq:1.8}), our method is iterative and in each iteration we partition the variables into
two part, namely the inactive set and the active set. Suppose the current iteration is  $l$, then we are going
to compute the $l+1$th iteration.

Based on the idea of fixed set in \cite{Kim1} and the theory in
\cite{Bertsekas}, we denote an inactive set as follows.
\begin{equation}\label{fixset}
I_+^l=\{i|0\leq h_i^l\leq\epsilon^l,\;[\nabla g(h)]_i>0\},
\end{equation}
where $\epsilon^l=\min(\epsilon,||h^l-\nabla g(h^l)||_F^2)$, $\epsilon$  is a small positive scalar. And $[\nabla g(h)]_i$
denotes the $i$th element of the gradient $\nabla g(h)$ of the object function $g(h)$.

For convenience, in some places we will slightly abuse notations, and say that $h_i^l\in I_+^l$ whenever $i\in I_+^k$.

Denote the active variables and the inactive variables at the
current iteration by $\xi^l$ and $\eta^l$, respectively. Assuming
that the $h^l$ and $\nabla g(h^l)$ are partitioned as follows:
\begin{equation}
h^l=\left[\begin{array}{ll}
\xi^l\\
\eta^l
\end{array}\right],
\end{equation}
where $\xi^l\notin I_+^l$ and $\eta^l\in I_+^l$. That is to say $\xi^l$ is the active part of the current iteration. Then
we compute the projection $\xi$ by the following equation:
\begin{equation}
\xi=P[\xi^l-\alpha\bar{D}^l\nabla g(\xi^l)],
\end{equation}
where $\alpha\geq0$ is the iteration step-size, $\bar{D}^l$ is a
gradient scaling matrix, $P[a]$ is the positive projection, i.e.,
\begin{equation}
P[a]=\left\{{\begin{array}{ll}
 a,\ \ \ \ a>0\\
0,\ \ \ \ a\leq0
 \end{array}}\right..
\end{equation}

Then we update the current result $h^l$ by the following rules
\begin{equation}\label{update}
h^{l+1}=\left[{\begin{array}{ll}
\xi\\
\eta^l
\end{array}}\right]=
\left[{\begin{array}{cc}
P[\xi^l-\alpha\bar{D}^l\nabla g(\xi^l)]\\
0
\end{array}}\right],
\end{equation}
where the right equation uses the fact that $\eta^l$ is fixed to
zeros, it can be comprehended as that this part satisfies the KKT
optimal condition Eq.(\ref{KKT}), then in the next iteration we need
not to update this part. After updating the $h^{l+1}$ we can compute
the $\nabla g(h^{l+1})$ and updating the active set $I_+^{l+1}$ to
obtain the next iterate result. In the whole iteration we will
adhere to the iteration rule until the result satisfies the stop
criterion.

\subsection{The Approximate to $\bar{D}^l$}

As the size of $\xi^l$ and $\eta^l$ changes at each iteration, the computation of the matrix $\bar{D}^l$ is not an
easy task. Note that the curvature information of $\xi^l$ is received from the curvature information of $h^l$.
Due to this fact, the gradient matrix $\bar{D}^l$ can be obtained by taking the proper sub-matrix of the
matrix $\bar{D}$ to avoid this task. Where $\bar{D}$ is a matrix that carry the curvature information of
the vector $h^l$. In realization of the method we can try to eliminate the curvature information of $\eta^l$.
Then the update rule (\ref{update}) can be regarded as follow:
\begin{equation}
h^{l+1}=P[h^l-\alpha\bar{D}\nabla g(h)].
\end{equation}

In \cite{Lin}, C. Lin used the BFGS update method to approximate the
Hessian matrix of the object function $g(h)$. The BFGS method is
well-established, and only uses the gradient information of the
object function. But BFGS is time consuming. Many researchers have
experimentally observed that the symmetric rank-one(SR1) rule
performs better than BFGS quasi-Newton update rules \cite{Oztoprak}.
In this paper we will use the SR1 rule to approximate the Hessian
matrix of $g(h)$ to improve the the converge speed.

Suppose $H^l$ is the current approximate of the Hessian matrix, then
using the SR1 to approximate the next Hessian matrix, we have that
\begin{equation}
H^{l+1}=H^l+\frac{(\omega^l-H^lv^l)(\omega-H^lv^l)^T}{(\omega-H^lv^l)^T\omega^l},
\end{equation}
where $\omega^l=\nabla g(h^{l+1})-\nabla g(h^l)$ and
$v^l=h^{l+1}-h^l$. Let $D^l$ be the inverse of the $H^l$, then it
can be achieved by using the Sherman-Morrison-Woodbury formula
\cite{Golub}. We can obtain that
\begin{equation}\label{inverse}
D^{l+1}=D^l+\frac{(v^l-D^l\omega^l)(v^l-D^l\omega^l)^T}{(v^l-D^l\omega^l)^T\omega^l}.
\end{equation}

\subsection{Line Search Strategy for Step-size}

From (\ref{update}), $\alpha>0$, we are considering to find $\alpha$
with the largest function reduction. For this purpose, we choose the
searching rule as follows:
\begin{equation}\label{stepsize}
\alpha=arg\min_{\alpha\geq 0}g(P[h^l-\alpha\nabla g(h^l)]).
\end{equation}

Thus, we obtain the following symmetric rank-one NNLS algorithm (see
Algorithm I).

\begin{tabular}{l}
\hline
Algorithm I (the SR1 algorithm for NNLS)\\ \hline
\textbf{Step 1}. \%start data\\
\ \ \ \ choose $h^0\geq 0,\; Denote\ the\ maximun\ number\ of\ iteration:maxiter,\;$\\
\ \ \ \ Set $l=0$; \\
\textbf{Step 2}. \\
\ \ \ \ compute the active set using the rule (\ref{fixset});\\
\textbf{Step 3}. \% update the current solution\\
\ \ \ \ compute the active set using Eq.(\ref{fixset}); \\
\ \ \ \ current solution by Eq.(\ref{update}), Eq.(\ref{inverse}) to approximate the inverse of  \\
\ \ \ \ Hessian matrix and Eq.(\ref{stepsize});\\
\textbf{Step 4}.\\
\ \ \ \ \textbf{if }the current result reach the stop criterion then stop, the current \\
\ \ \ \ solution is the final solution;\\
\ \ \ \ \textbf{else} goto Step 2 and Step 3.\\
\hline
\end{tabular}
\setcounter{equation}{0}
\renewcommand{\theequation}{3.\arabic{equation}}

\section{A Symmetric Rank-one Quasi Newton Method for NMF}

In the former part, we have discussed the algorithm for NNLS, then
in this part we are going to talk about the algorithm for NMF. Note
that, in Section 1, we have analysed that the NMF problem can be
resolved into several non-negative least squares. Due to this
property, we can apply the symmetric rank-one method to our NMF
problem.

\subsection{Applying the SR1 Directly to NNLS Subproblems}

For Eq.(\ref{Eq:1.5}), we can do vectorization of the matrix $H$,
and obtain the {\small\begin{eqnarray}\label{vec}
g(H)&=&\frac{1}{2}||V-WH||_F^2\nonumber\\
&=&\small{ \frac{1}{2}vec(H)^T\left({\begin{array}{*{21}c}
{W^TW} & {} & {}\\
{} & {\ddots} &{}\\
{} & {} & {W^TW}
\end{array}}\right)vec(H)-tr((W^TV)^TH)+\frac{1}{2}tr(V^TV)},
\end{eqnarray}}
where $tr(B)$ is the trace of matrix B. Then using the quasi-Newton idea, the iteration of $vec(H)$ is
\begin{equation}\label{iter}
vec(H^{l+1})=P[vec(H^l)-\alpha^l\bar{D}\nabla g(H)],
\end{equation}
where $\alpha^l$ is the step-size of the current iteration, $\bar{D}$ is the gradient scaling matrix and $\nabla g(H)$
is the gradient of (\ref{vec}).

As we all known, Eq.(\ref{vec}) is actually a non-negative least
square problem, then we can apply the method discussed in the former
part to this problem directly. First let the inactive set of
(\ref{vec}) be defined as
\begin{equation}
I^l_H=\{i|0\leq (vec(H))_i\leq\epsilon^l,\nabla g(H)_i>0\},
\end{equation}
where $\epsilon^l$ can be obtained by the same method of the former part.
In Eq.(\ref{iter}), the gradient scaling matrix $\bar{D^l}$ is approximated by the symmetric rank-one quasi Newton method.
The step-size $\alpha^l$ is done by the line search. For Eq.(\ref{Eq:1.4}), we can do by the same method.
The object function corresponding to $W$ is
{\small\begin{eqnarray}\label{vec2}
g(W)&=&\frac{1}{2}||V^T-H^TW^T||_F^2\nonumber\\
&=&\frac{1}{2}vec(H)^T\left({\begin{array}{*{21}c}
{HH^T} & {} & {}\\
{} & {\ddots} &{}\\
{} & {} & {HH^T}
\end{array}}\right)vec(W^T)-tr((HV^T)^TW^T)+\frac{1}{2}tr(V^TV).
\end{eqnarray}}

\subsection{The Algorithm for NMF}

Next, based on the Algorithm I, we give the symmetric Rank-one
quasi-Newton algorithm for NMF as follows.

\begin{tabular}{l}
\hline Algorithm II( the  algorithm for NMF)\\ \hline
Initialization\\
\ \ \ \ input matrix $V$, the maximum number of iteration: maxiter;\\
\ \ \ \ initial matrix $W^0$\\
\ \ \ \ the approximate rank $K,\ 1\leq K\leq\min{m,n}$\\
for $i=1:maxiter$\\
\ \ \ \ (1) $W=W^l,\ H^{old}=H^l$\\
\ \ \ \ (2) compute the matrix $H^{new}$ by using the Algorithm 1, where the object function\\
\ \ \ \ \ \ \ \ \ is Eq.(\ref{vec})\\
\ \ \ \ (3) $H^{old}=H^{new},H^{l+1}=H^{old}$\\
\ \ \ \ (4) $H=H^{l+1},W^{old}=W^l$\\
\ \ \ \ \ \ \ \ \ compute matrix $W^{new}$ by using Algorithm 1, where the object function\\
\ \ \ \ \ \ \ \ \ is Eq.(\ref{vec2})\\
\ \ \ \ (5) $W^{old}=W^{new},W^{l+1}=W^{old}$\\
\ \ \ \ (6) if the stopping criteria is met, break\\
end for\\
\hline
\end{tabular}

The computation of the sub-problem (\ref{Eq:1.4}) and (\ref{Eq:1.5}) is hight cost. Each sub-problem requires an
iterative procedure, which can be regarded as an sub-iteration, and in our algorithm, there exist two
computation of the gradient:
\begin{equation}
{\begin{array}{ll}
\nabla f_H(W,H)=W^TWH-W^HV,\\
\nabla f_W(W,H)=WHH^T-VH^T.
\end{array}}
\end{equation}

During to computation, we can compute $W^TW, W^TV$ and $H^TH,VH^T$ in the outer iteration.
\setcounter{equation}{0}
\renewcommand{\theequation}{4.\arabic{equation}}

\section{Numerical Experiments}

In this part, we present some experiment results from both rand data
and real world data and compare our method to the other six
following algorithms.

\begin{enumerate}
  \item \textbf{fnmae}: the projected quasi-Newton method in \cite{Kim1}.
  \item \textbf{pnm}: alternating nonnegative least squares using the projected Newton method in \cite{gong}.
  \item \textbf{AS}: the active set method in \cite{Jkim2}.
  \item \textbf{alsq}: the alternating nonnegative least squares method in \cite{Berry}.
  \item \textbf{nnmalq}: the alternating nonnegative least squares method in \cite{Lee1}.
  \item \textbf{nmf}: the projected gradient method in \cite{Lin}.
  \item \textbf{SR1}: the symmetric rank-one quasi Newton method in this paper.
\end{enumerate}

In our experiments, we initialize all the method randomly, and show plots of the relative error against the number of iteration or the
time of iteration, where the relative error of approximation is $||V-WH||_F/||V||_F$.

\subsection{Synthetic Data Experiment }

We generate the synthetic data  randomly, in our method we test the randomly generated matrix $V\in\mathbb{R}^{n\times m}\
(n=200,\ m=40)$, the approximate rank is $k=10$ and another matrix $V\in\mathbb{R}^{n\times m}\ (n=2000,\ m=800)$,
which we set the approximate rank
$k=5,10,20$, respectively. In this experiment we all done ten times of the random matrix, then we take the average result of the
ten experiment data.

\begin{figure}[htbp]\label{fig1}
\centerline{\includegraphics[width=3.64in,height=2.88in]{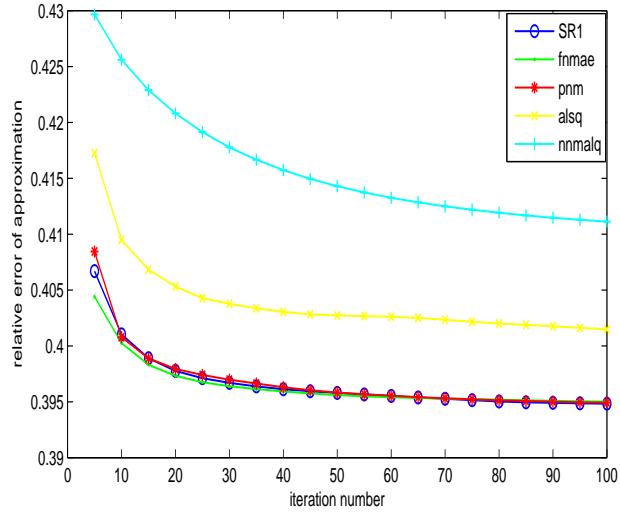}}
\caption{\small Relative error of approximate against iteration number count for SR1, fnmae, pnm,alsq,nnmalq with the
matrix$V\in\mathbb{R}^{200\times 40}$ for the approximate rank $k=10$.}
\end{figure}

\begin{figure}[htbp]\label{fig2}
\centerline{\includegraphics[width=6.64in,height=3.50in]{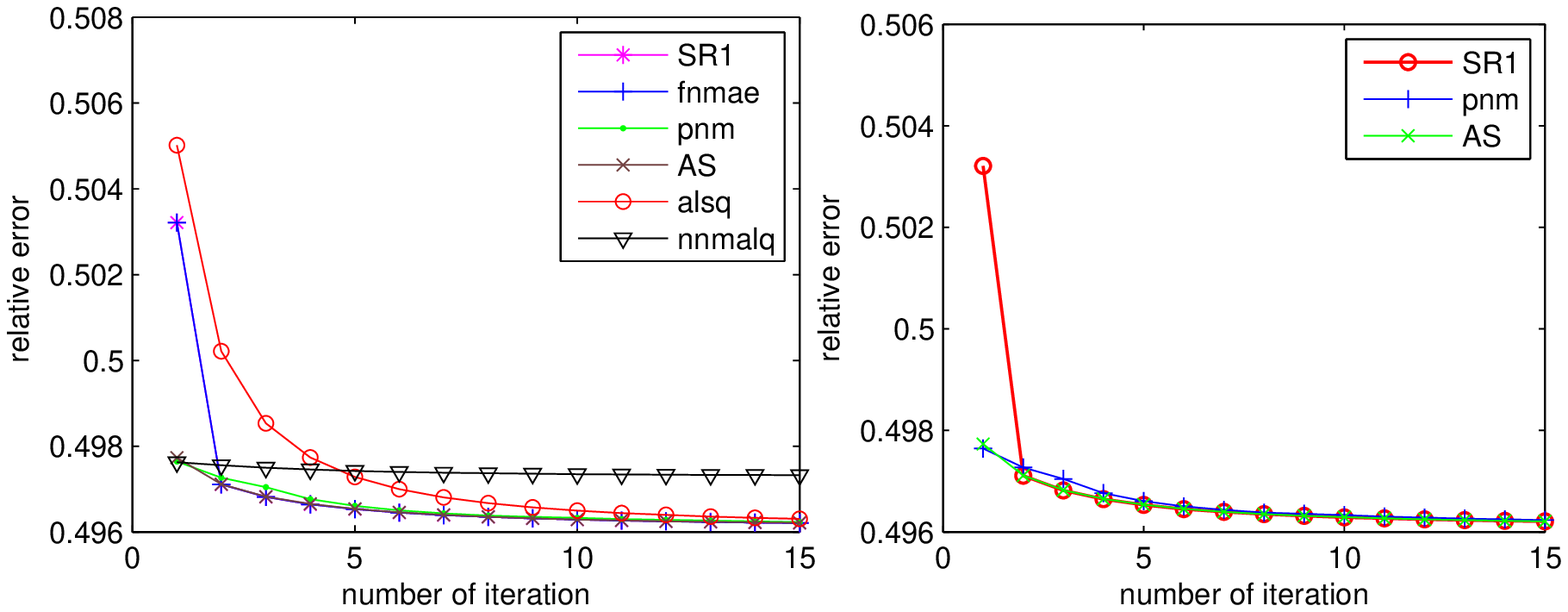}}
\caption{\small Relative error of approximate against iteration number count for SR1, fnmae, pnm,AS,alsq,nnmalq with the
matrix $V\in\mathbb{R}^{2000\times 800}$ for the approximate rank $k=5$.}
\end{figure}

\begin{figure}[htbp]\label{fig3}
\centerline{\includegraphics[width=6.64in,height=3.50in]{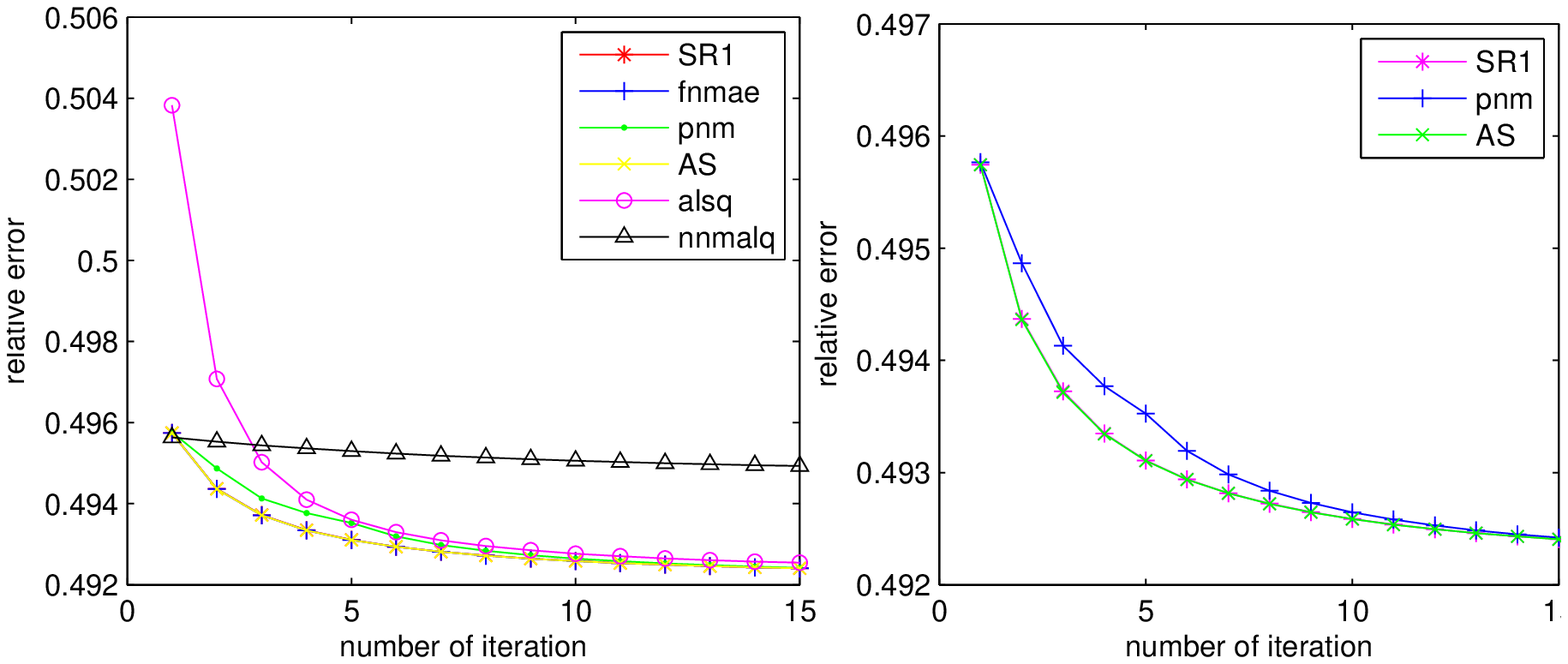}}
\caption{\small Relative error of approximate against iteration number count for SR1, fnmae, pnm,AS,alsq,nnmalq with the
matrix $V\in\mathbb{R}^{2000\times 800}$ for the approximate rank $k=10$.}
\end{figure}

\begin{figure}[htbp]\label{fig4}
\centerline{\includegraphics[width=6.64in,height=3.50in]{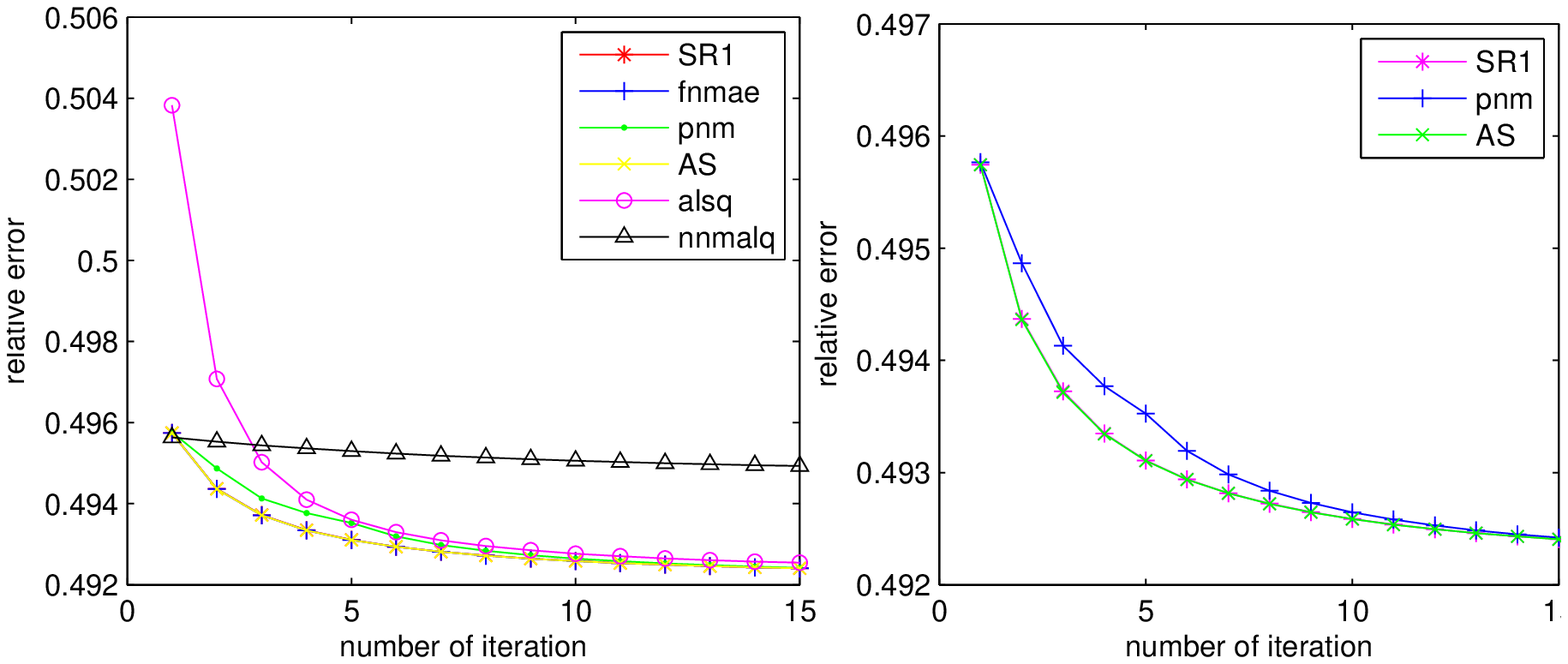}}
\caption{\small Relative error of approximate against iteration number count for SR1, fnmae, pnm,AS,alsq,nnmalq with the
matrix $V\in\mathbb{R}^{2000\times 800}$ for the approximate rank $k=20$.}
\end{figure}

\begin{figure}[htbp]\label{fig5}
\centerline{\includegraphics[width=6.64in,height=3.00in]{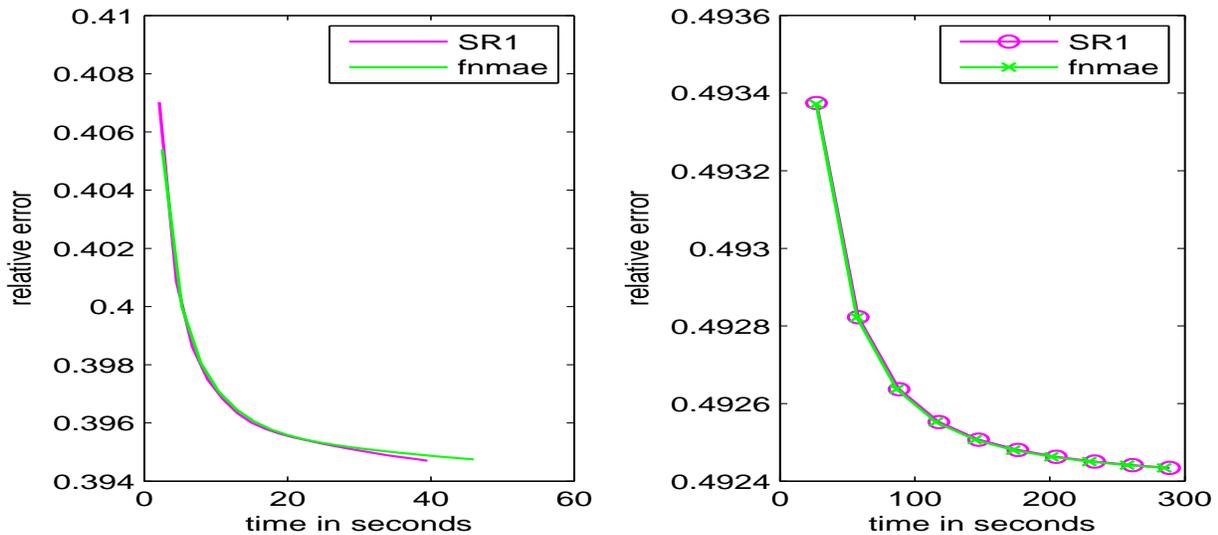}}
\caption{\small Relative error of approximate against time count for SR1, fnmae with the
matrix $V\in\mathbb{R}^{200\times 40}$ in the left hand side and $V\in\mathbb{R}^{2000\times 800}$
 in the right hand side for the approximate rank $k=10$.}
\end{figure}

\begin{figure}[htbp]\label{fig5}
\centerline{\includegraphics[width=4.00in,height=2.40in]{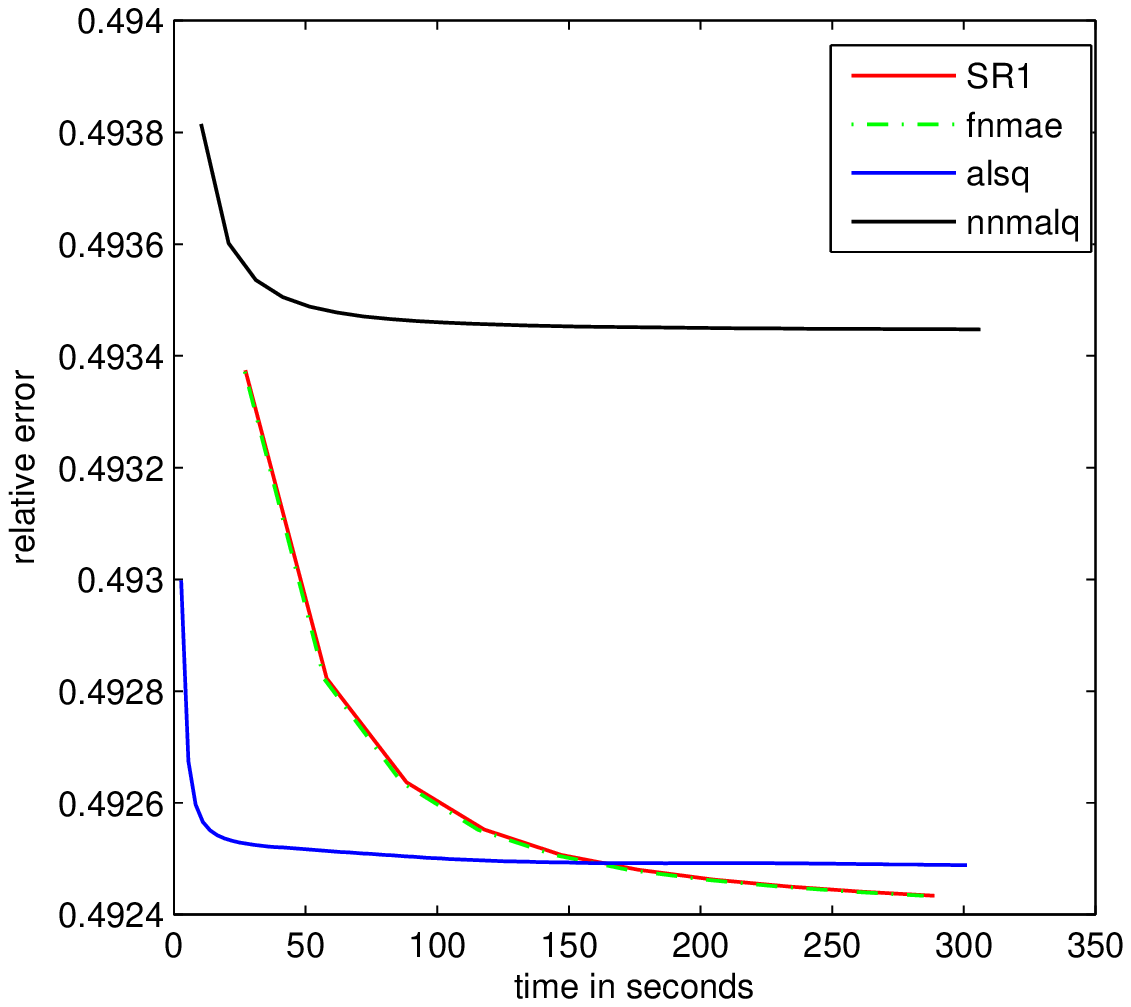}}
\caption{\small Relative error of approximate against time count for SR1, fnmae, alsq, nnmalq with the
matrix  $V\in\mathbb{R}^{2000\times 800}$ for the approximate rank $k=10$.}
\end{figure}

Figure 1 explains the random matrix $V\in\mathbb{R}^{200\times 40}$
for the rank $k=10$, and from the figure we can learn that,
comparing with fnmae \cite{Kim1}, the object function decreases more
per iteration, and our method costs less time to reach the same
relative error. As a whole the relative errors of approximate of the
 SR1, fnmae and pnm are very similar.

Figure 2-4 are the numerical results of the random matrix
$V\in\mathbb{R}^{2000\times 800}$ and setting the approximate rank
$k=5,10,20$, respectively. The figures in the left side is the
comparison of SR1, fnmae, pnm, AS,alsq and nnmalq, to have a clear
look at  similar numerical results of SR1, pnm and AS methods we
present the comparison of the three methods in the right hand side.
From these figures we can learn that, with the random initialization
the relative error of approximate of SR1, fnmae, pnm and AS are
similar, but our method has a slight advantage.

Figure 5 and 6 are the plots of relative error against the running
time, we can learn that our method converges faster than the fnmae
method, and much faster than the nnmalq and alsq methods.

The numerical results show that our symmetric rank-one quasi-Newton method improves the efficiency of the quasi-Newton
method every iteration, and cuts down the time of each iteration. Comparing with other methods, our SR1 method decreases
faster in each iteration. And the nnmalq method decreases  slow in every iteration, so become less competitive. Since
 in our experiment the error of nmf method is much lager than the other six methods mentioned above, we did not
plot the error curve of this method.

\subsection{ Application to imaging processing}

NMF was originally motivated by Lee and Seung \cite{Lee2} using an image processing application. Many others have
also considered NMF to image processing, face recognition application, model recognition application and signal processing
application. In this part we are going to do some numerical experiments of image processing application.
Our experiment is done on four random choosen faces\footnote{http://www.cl.cam.ac.uk/research/dtg/attarchive/facedatabase.html}
 and the size of each face image is $92\times 112$ with 256
gray per pixel. We give out the the reconstructed image of running 10, 20, 50
iterations, respectively. In this experiment, we set the approximate rank $k=14$, and the images are
taken randomly from the 5 random initial matrix in each each iteration. The first image in each row is the original image
, followed by reconstructions obtained by SR1, fnmae, pnm, AS, alsq, nnmalq and nmf.

From figures 7 and 8 , we can learn that, after 10 and 20 iterations
SR1 can reconstruct better than fnmae, AS, pnm procedures, and much
better than alsq, nnmalq and nmf procedures. Figure 9 illustrates
that, the images obtained via SR1, AS, fnmae, pnm methods have
similar quality, and the image obtained by alsq is still vague.

\begin{figure}[htbp]\label{fig5}
\centerline{\includegraphics[width=6.54in,height=2.80in]{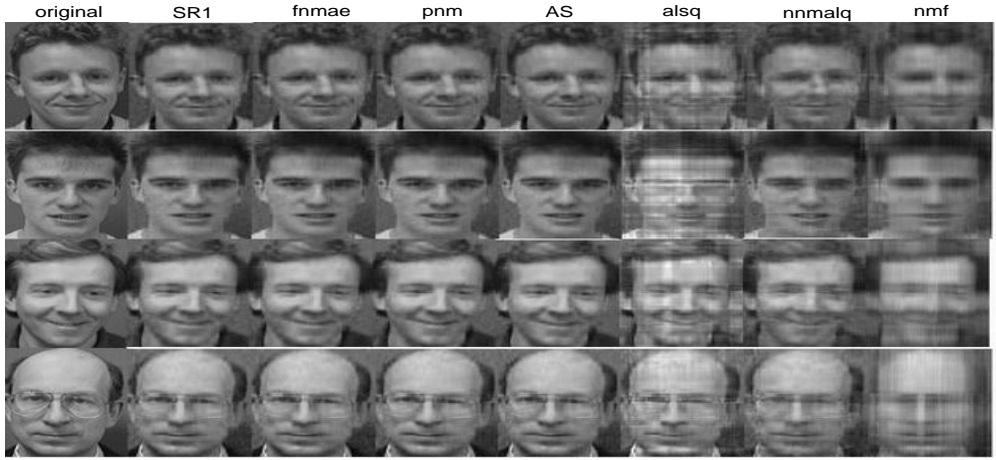}}
\caption{\small Reconstructed images after 10 iterations.}
\end{figure}

\begin{figure}[htbp]\label{fig6}
\centerline{\includegraphics[width=6.54in,height=3.40in]{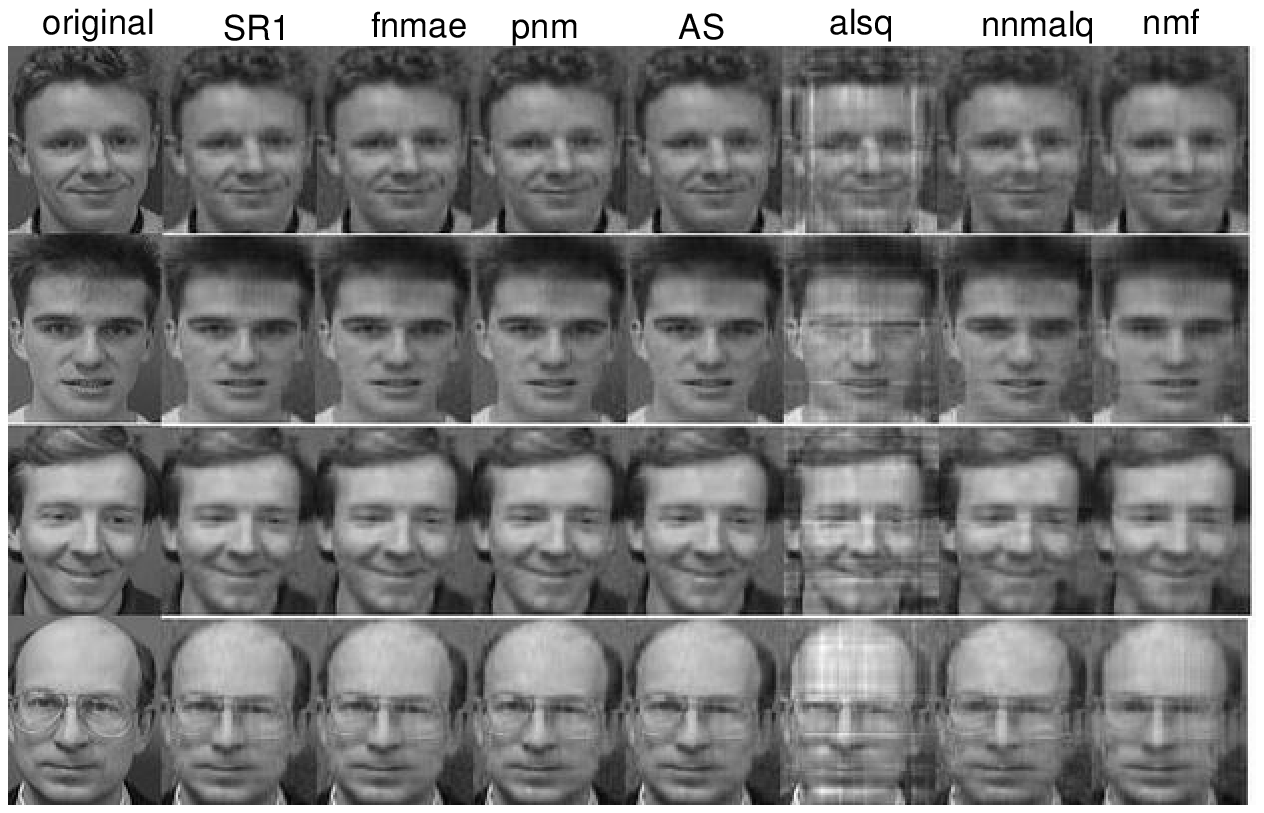}}
\caption{\small Reconstructed images after 20 iterations}
\end{figure}

\begin{figure}[htbp]\label{fig7}
\centerline{\includegraphics[width=6.54in,height=3.40in]{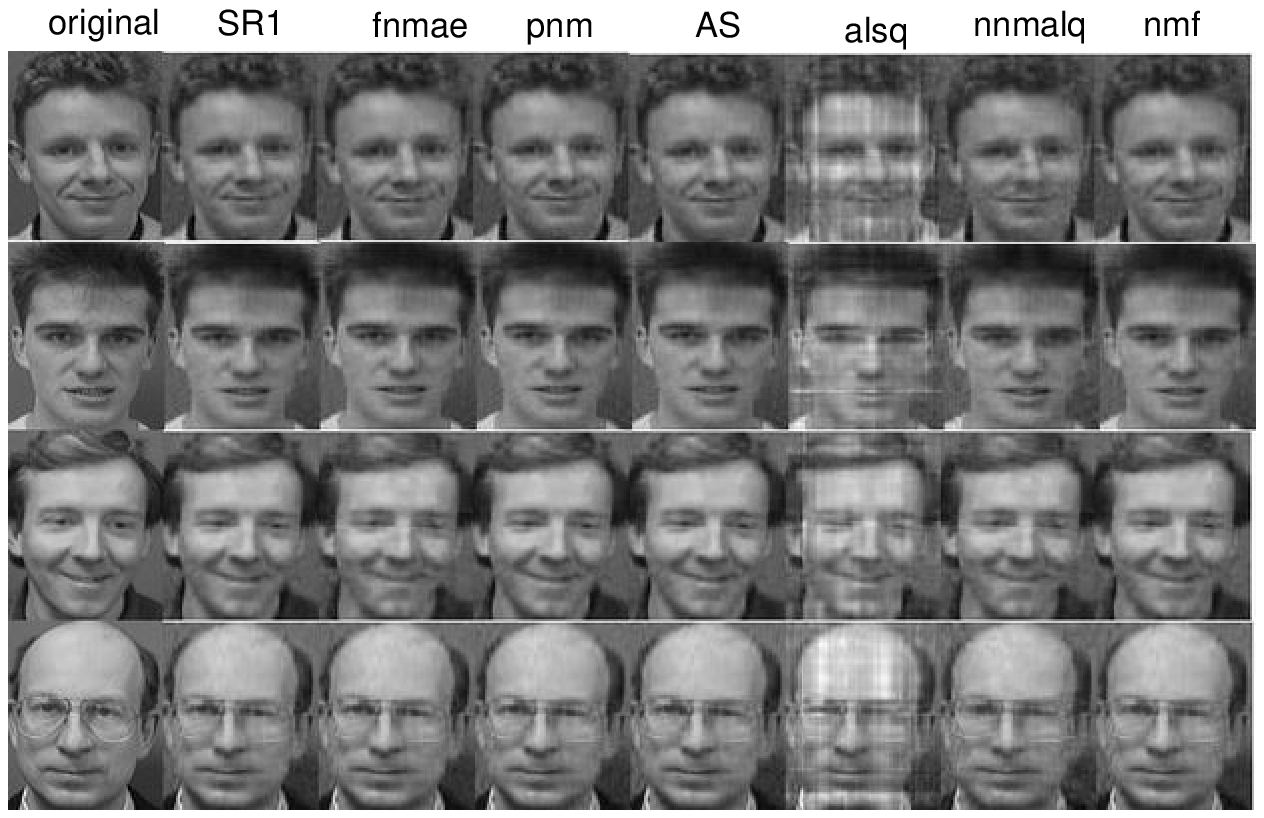}}
\caption{\small Reconstructed images after 50 iterations}
\end{figure}

\subsection{ Application to Text Clustering}

In order to test the application to text analysis, we are going to apply our method to text
clustering, and compare with the other methods. We show the numerical results of the above methods
except the nmf method on four text datasets, which are high dimensional and sparse. The text data
is collected from newspapers, such us LA times, San Jose Mercury and so on\footnote{http://www.shi-zhong.com/software/docdata.zip}.
In the numerical experiments
we all set the approximate rank $k$ to the number of classes, and we all initialize the matrix
randomly. As the same with the random data, we present the relative error of approximate
against the number of iterations.

\begin{figure}[htbp]\label{fig8}
\centerline{\includegraphics[width=5.04in,height=3.30in]{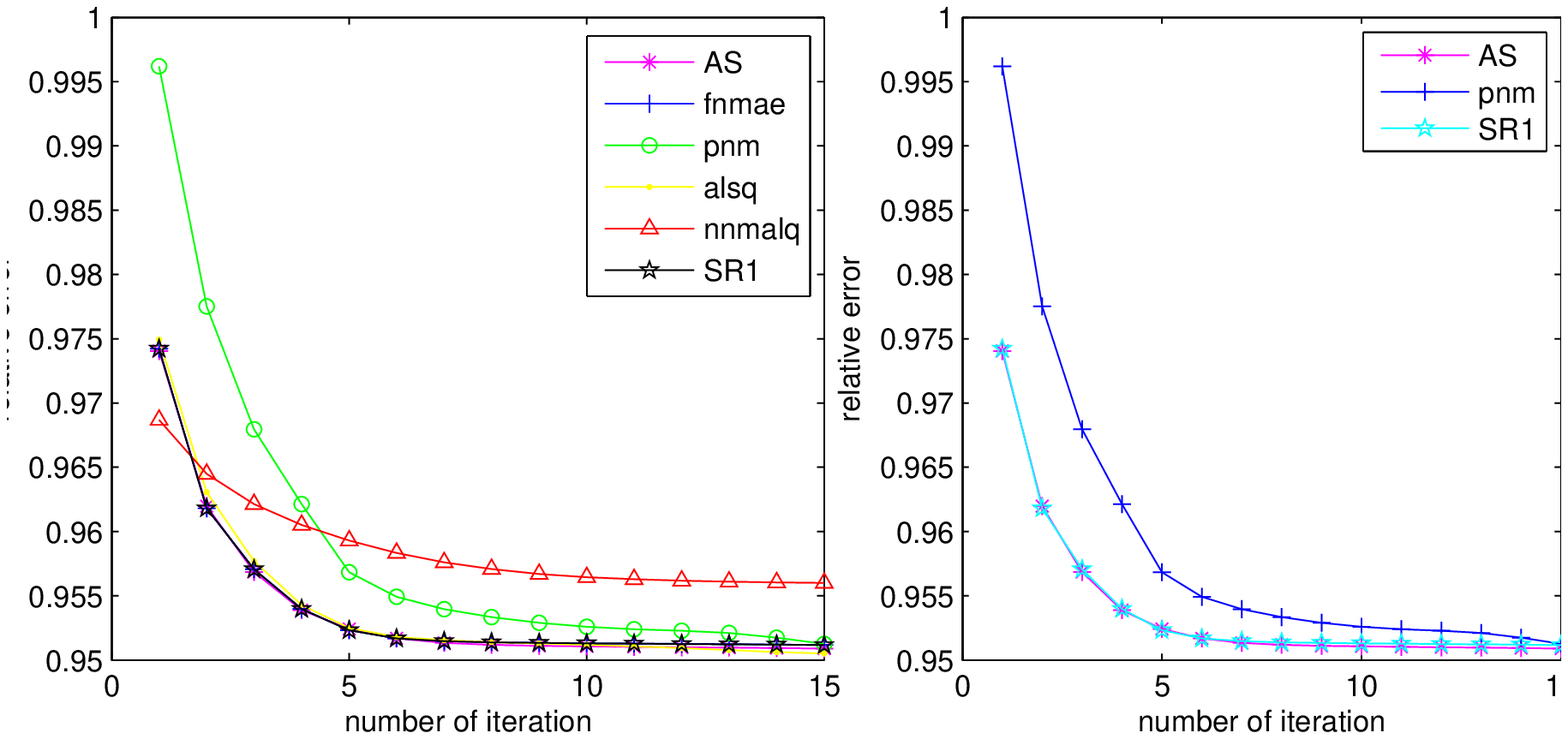}}
\caption{\small Relative error of approximate against the number of iteration for the text data classic 300}
\end{figure}

\begin{figure}[htbp]\label{fig9}
\centerline{\includegraphics[width=5.54in,height=2.50in]{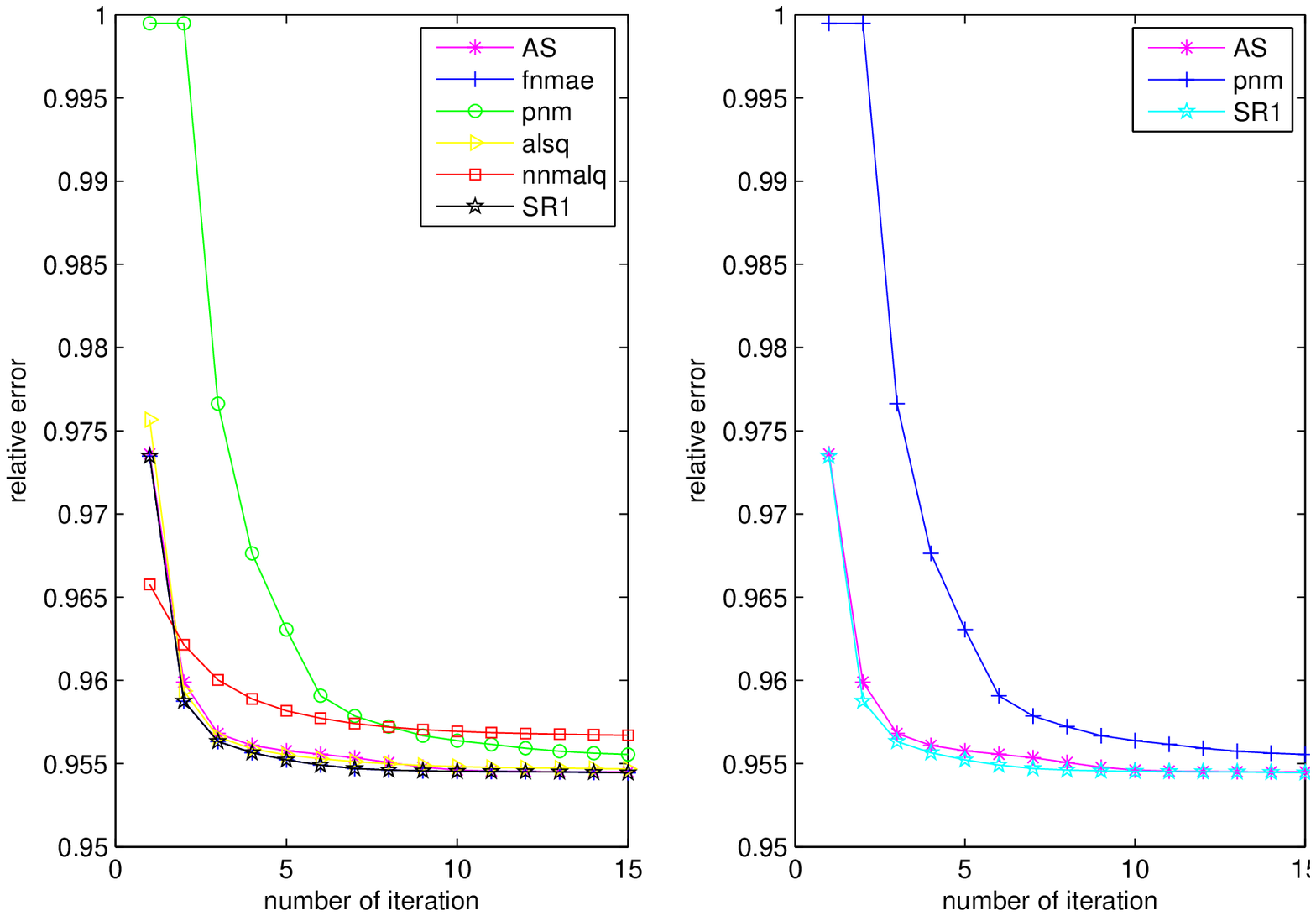}}
\caption{\small Relative error of approximate against the number of iteration for the text data classic 3891}
\end{figure}

\begin{figure}[htbp]\label{fig10}
\centerline{\includegraphics[width=5.04in,height=3.30in]{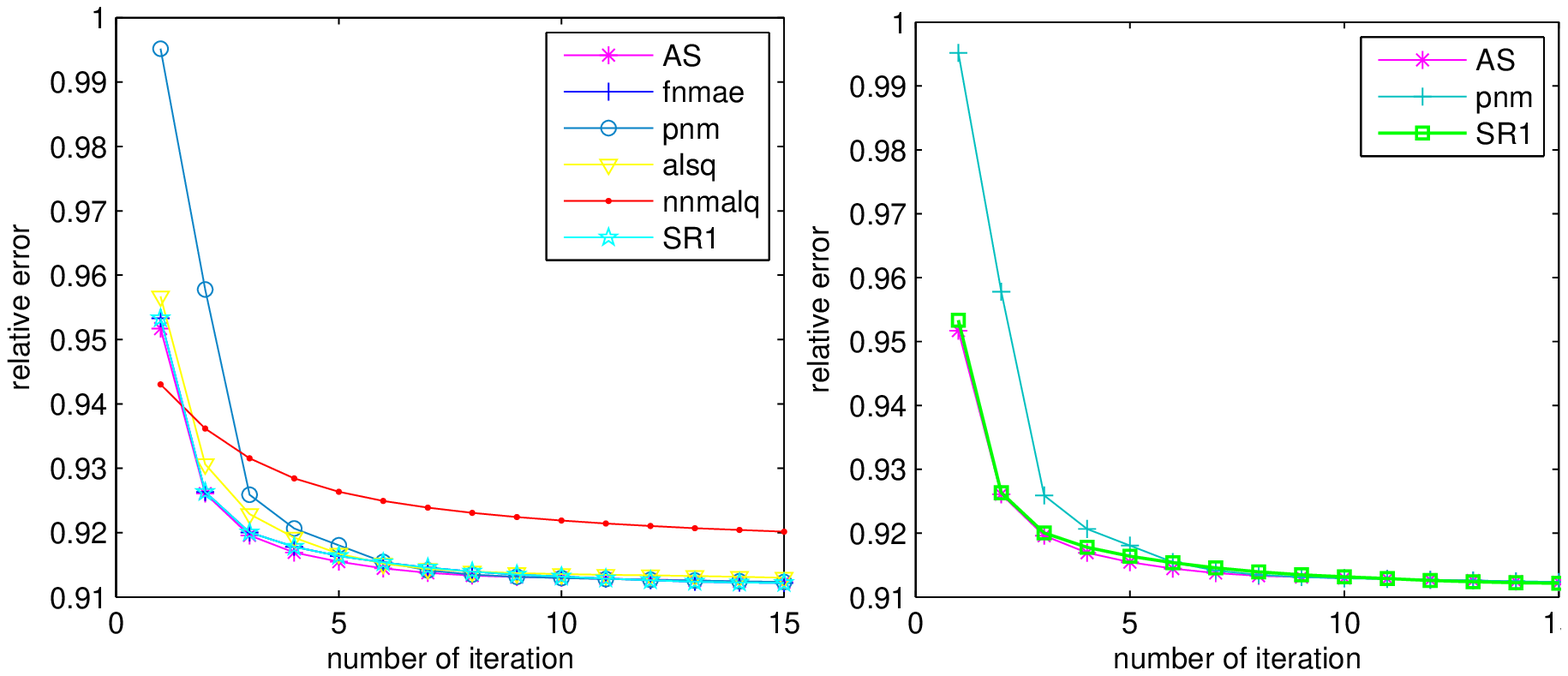}}
\caption{\small Relative error of approximate against the number of iteration for the text data hitech}
\end{figure}

\begin{figure}[htbp]\label{fig11}
\centerline{\includegraphics[width=5.04in,height=3.30in]{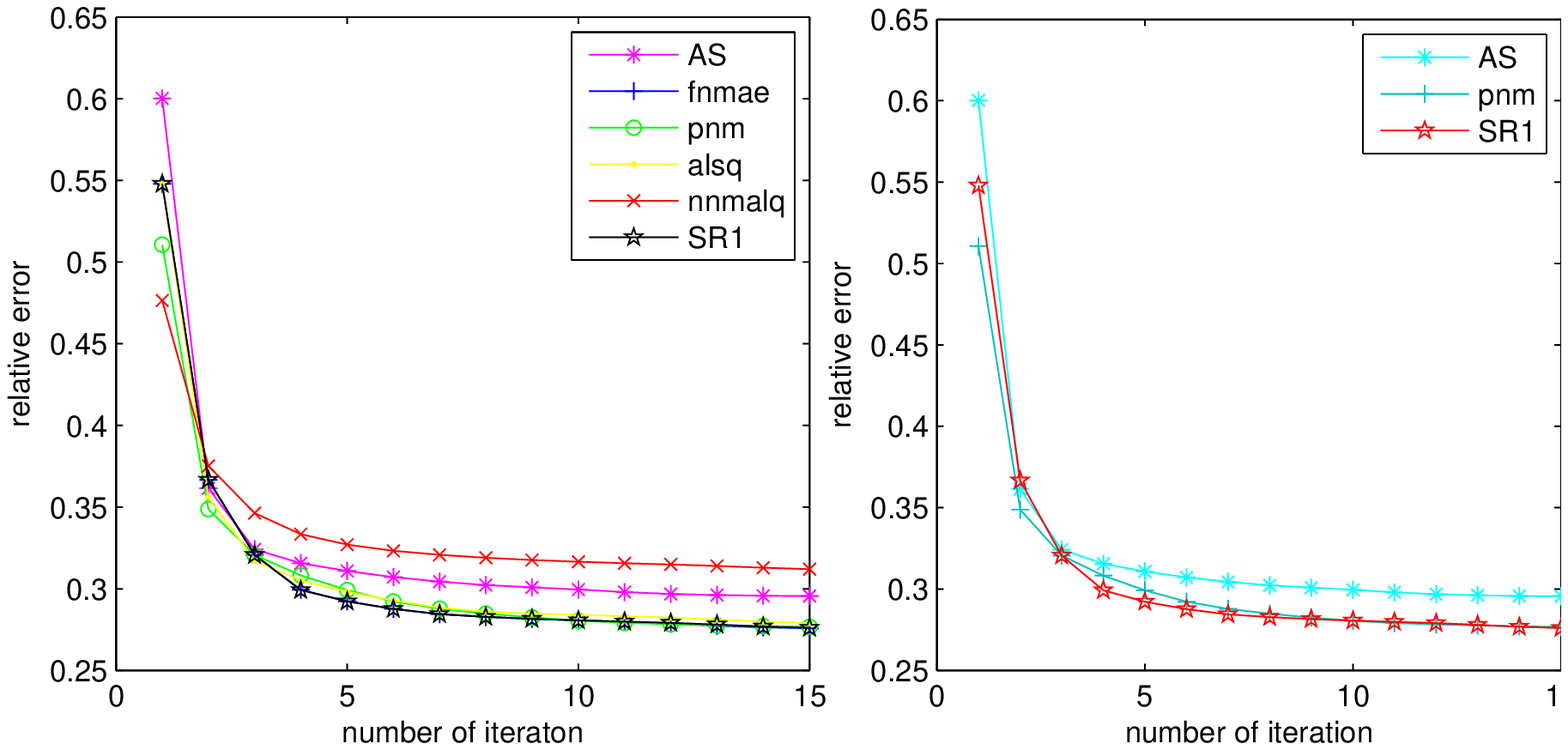}}
\caption{\small Relative error of approximate against the number of iteration for the text data tr45}
\end{figure}

We present the comparison of the six methods mentioned above except
the nmf method in the left hand plot. The right hand plot is to
outstand the comparison of the SR1, AS and pnm methods. From figures
10-13, we can known that with the same random initial matrix, at the
beginning several iterations our method is much better than pnm, and
decrease much faster than the other methods. After 5 iterations, the
quality of the methods are more and more similar.

\setcounter{equation}{0}
\renewcommand{\theequation}{5.\arabic{equation}}

\section{Conclusion}

In this paper we present an algorithm for NMF, it is different from the fnmae method \cite{Kim1} in the following
aspects:

 (1) The active set in our method is a relax form, while the active set in \cite{Kim1} is an hard form. The fact that
 an active set of hard for exhibits undesirable discontinuity at the boundary of the constraint set has been showed
 in Bertsekas\cite{Bertsekas}. This is harm to the convergence rate, so we use the relax form to
 avoid this problem.

 (2) In addition, the fnmae method \cite{Kim1} using the BFGS method to approximate the Hessian matrix, we use the symmetric
 rank one method to approximate the hessian matrix. The symmetric rank one method has been experimentally shown
 better than the BFGS method.

 What is more, our method is also different from the pnm method in \cite{gong} for the aspect that the approach of the approximating the
 Hessian matrix. In \cite{gong}, using the symmetric property of the Hessian matrix approximates the inverse of the Hessian matrix by
 the Choleskey factorization, while in our method we approximate the
inverse of the Hessian matrix by the symmetric rank-one method.
Numerical experiments show that the object function decreases more
per iteration in our method than in the pnm method \cite{gong} in
some cases.

All in all, we propose the symmetric rank one quasi Newton method
for the NMF. This method maintains the decrease speed per iteration
and decrease the computational time. From the experimental results
both the synthetic data and real world data, our method performs
very well. Due to the wide use of the non-negative matrix
factorization, to learn more effective methods is necessary. In
addition, there are many factors which influence the efficiency of
the algorithm for NMF, such as the approximate rank, sparseness and
so on. So researching more efficiency algorithm for NMF by
approximating rank and the stop criterion is necessary in the
future.

\textbf{Acknowledgements} \emph{We sincerely thank the authors (see
\cite{Lin}, \cite{Berry}, \cite{Kim1}, \cite{Jkim2}, \cite{Lee1},
\cite{gong}) of the above six methods for their experiments codes.
With the help of them, it is more convenient for us to do these
numerical experiments.}

\end{document}